\documentclass[12pt]{article}
\usepackage{amsfonts}
\usepackage{amssymb}
\usepackage{mathrsfs}
  \date{ }
   \title{{\bf A kind of infinite-dimensional Novikov algebras and its realization}
  \thanks{Supported by  NNSF of China (No.10871057),  Natural
Science Foundation of  Jilin province (No.201115006), Scientific
Research Foundation for Returned Scholars
    Ministry of Education of China and the Fundamental Research Funds for the Central Universities.\quad \quad\quad\quad\quad
    \quad\quad\quad \quad\quad Corresponding author(L. Chen):
     chenly640@nenu.edu.cn}}
   \author{{Liangyun Chen, Yao Ma, Haijun Yu}
   \\  Department of Mathematics, Northeast Normal University,
     \\ Changchun 130024, CHINA }
  
\begin{document}
 \baselineskip=0.3in
  \maketitle
 \begin{center}{\bf Abstract}\end{center}

  In this paper,  we  construct a  kind of infinite-dimensional Novikov algebras and
   give its realization by hyperbolic sine functions and  hyperbolic cosine functions.



\begin{center}
 {\bf{\S 1 \quad Introduction}}
\end{center}

 \noindent   The Hamilton operator is an important operator of the
 calculus of variations.  When I. M. Gel'fand and I. Ya. Dorfman [5-7] studied the following operator:
     $$H_{ij}=\sum_{k}c_{ijk}u_{k}^{(1)}+d_{ijk}u_{k}^{(0)}\frac{\rm d}{{\rm d}x}, \quad c_{ijk}\in{ \bf {C}},
   \quad d_{ijk}=c_{ijk}+c_{jik},$$
  they  gave the definition of Novikov algebras. Concretely,
  let $c_{ijk}$  be the structural coefficients,
   a product of $L=L(e_{0},e_{1}, \cdots )$ be  $\circ$ such that
  $$e_{i}\circ{e_{j}}=\sum{c_{ijk}e_{k}}.$$  Then the product is Hamilton operator if and only if
  $\circ$ satisfies£º
  $$(a\circ b)\circ c=(a\circ c)\circ b$$
  $$(a\circ b)\circ c+c\circ (a\circ b)=(c\circ b)\circ a+a\circ (c\circ b).$$

In 1987, E. I. Zel'manov[15] began to study Novikov algebras and
 proved that the dimension of finite simple  Novikov algebras over a
field of  characteristic zero is one. In algebras, what are paid
attention to by mathematician are  classifications and structures, but so far we haven't got the systematic theory for
general Novikov algebras. In 1992, J. M. Osborn [9-10] had  finished the
classification of  infinite simple  Novikov algebras with
nilpotent elements over  a field of  characteristic zero and finite
simple Novikov algebras with nilpotent elements over a field of
 characteristic $p>0$. In 1995, X. P. Xu [11-14] developed his theory and
 got  the classification of simple Novikov algebras
 over an algebraically closed field of  characteristic zero.
 C. M. Bai and  D. J. Meng [1-3]  has serial work
 on low dimensional  Novikov algebras, such as the  structure  and
 classification. We construct  two kinds of  Novikov algebras [4]. Recently,  people obtain some  properties in  Novikov superalgebras [8, 16].

     In this paper,  we  construct  a new infinite-dimensional Novikov algebras  and  give its realization by hyperbolic sine functions and hyperbolic cosine functions.

 \noindent {\bf Definition 1.1 }\quad {\it Let $({\cal A}, \circ)$ be an algebra over $\bf F$
such that:
$$ (a, b,c)=(b, a, c),\eqno(1.1)$$
$$ (a\circ b)\circ c=(a\circ c)\circ b, \quad
 \forall a, b,c\in{\cal A}, \eqno(1.2)$$ then $\cal A$ is called a  Novikov algebra  over $\bf F$.}

 \noindent{\bf Remark 1.2 }\quad  (1)\, {\it Condition $(1.1)$ is  usually written by
  $$a\circ(b\circ c)-(a\circ b)\circ c=b\circ(a\circ c)-(b\circ a)\circ c. \eqno(1.3)$$ }
  (2)\quad {\it  An algebra $\cal A$ is called {a left
    symmetric algebra} if it only satisfies $(1.1)$. It is clear that left symmetric algebras  contain  Novikov algebras.}

  \noindent{\bf Remark 1.3 }\quad  (1)\, {\it  If $({\cal A}, \circ)$ is a  left symmetric
 algebra satisfying
 $$[a, b]=a\circ b-b\circ a,\quad \forall a, b\in{\cal A},\eqno(1.4)$$
    then $({\cal A}, [,])$ is a Lie algebra. Usually, it is called  an adjoining  Lie algebra. }

  (2)\quad {\it Let $({\cal A}, \cdot)$  be a commutative algebra, then $({\cal A},d_0, \circ)$ is a
   Novikov algebra if
  $d_0$ is a derivation of $\cal A$ with a bilinear  operator
 $\circ$  such that
 $$a\circ b=a\cdot d_0(b), \quad \forall a, b\in{\cal A}.\eqno(1.5)$$ }
 \begin{center}{\bf \S 2 \quad  Main results}\end{center}

 \noindent{\bf Lemma 2.1}\quad {\it Let \{ $b_0$, $a_1$, $b_1$, $a_2$, $b_2$, $\cdots$
 $a_n$,  $b_n$, $\cdots$ \} be  a basis of the linear space $\cal A$
 over a field  $\bf F$  of characteristic $p\neq 2$ satisfying
$$\left\{\begin{array}{l}a_ma_n=\displaystyle\frac{1}{2}(b_{m+n}-b_{m-n}),\cr \cr
 b_mb_n=\displaystyle\frac{1}{2}(b_{m+n}+b_{m-n}),\cr  \cr
 a_mb_n=b_na_m=\displaystyle\frac{1}{2}(a_{m+n}+a_{m-n}),\cr\end{array}\right.\eqno(2.1)$$
 where $b_{-m}=b_m$, $a_{-m}=-a_m$.  Then $\cal A$ is a commutative and associative algebra.}

  {\it Proof}.  It is clear that  $\cal A$ is a  commutative algebra over  ${\bf F}$.
    \begin{eqnarray*}&&(a_k, a_n, a_m)=a_k(a_na_m)-(a_ka_n)a_m\\
&=&\frac{1}{2}a_k(b_{m+n}-b_{n-m})-\frac{1}{2}(b_{k+n}-b_{k-n})a_m\\
&=&\frac{1}{4}(a_{k+m+n}+a_{k-m-n}-a_{k+n-m}-a_{k-n+m}\\
&&-a_{m+k+n}-a_{m-k-n}+a_{m+k-n}+a_{m-k+n})\\
&=&0.\end{eqnarray*} Similarly, we have that
 $(b_k, b_n, b_m)= (a_k, a_n, b_m)= (a_k, b_n, a_m)= (b_k, a_n, a_m)= (b_k, b_n, a_m)=
 (b_k, a_n, b_m)= (a_k, b_n, b_m)=0$. Then
$ (a, b, c)=0,  \forall a, b, c\in{\cal A}.$
The result follows.  \hfill$\Box$

 \noindent{\bf Corollary 2.2}\quad
{\it  $b_0$ of  Lemma 2.1 is a unity of $\cal A$.}

\noindent{\bf Lemma 2.3}\quad
 {\it Let $\cal A$ be a commutative and associative
algebra  satisfying  Lemma 2.1.  Then the following statements hold:

$1)$\quad If $D_0$ is a linear transformation of $\cal A$ such that
$$\left\{\begin{array}{ll}D_0(a_n)=nb_n,&n=1, 2, 3,\cdots,\cr
D_0(b_n)=na_n, &n=0, 1, 2, \cdots,\end{array}\right.\eqno(2.2)$$
then $D_0$ is  a derivation of $\cal A$.

$2)$\quad  If $aD_0$ is a linear transformation of $\cal A$ such that
$$
(aD_0)(b)=aD_0(b),  \forall a, b\in{\cal A},\eqno(2.3)$$
then  $aD_0$ is  a derivation of $\cal A$.

$3)$\quad ${\cal D}_1=\{aD_0|a\in {\cal A}\}$ is a subalgebra of Lie
algebra ${\rm Der}{\cal A}$.}

 {\it Proof}. \quad (1)\, We have
 $$D_0(a_na_m)=D_0\left(\frac{1}{2}(b_{n+m}-b_{n-m})\right)=\frac{1}{2}((m+n)a_{n+m}-(n-m)a_{n-m})$$
and
\begin{eqnarray*}
D_0(a_n)a_m+a_nD_0(a_m)&=&nb_na_m+ma_nb_m\\
&=&\frac{n}{2}(a_{n+m}-a_{n-m})+\frac{m}{2}(a_{n+m}-a_{m-n})\\
&=&\frac{1}{2}((m+n)a_{m+n}-(n-m)a_{n-m}).\end{eqnarray*}
   So $D_0$ is  a derivation of $\cal A$.

2)\quad For $\forall a, b, c\in{\cal A}$, we have
$$(aD_0)(bc)=aD_0(bc)=aD_0(b)c+abD_0(c)=(aD_0)(b)c+b(aD_0)(c),$$
so $aD_0$ is a derivation of $\cal A$.

3)\quad For $\forall a, b, c\in{\cal A}$, we have
\begin{eqnarray*}[aD_0, bD_0](c)&=&(aD_0)(bD_0)(c)-(bD_0)(aD_0)(c)\\
&=&aD_0(b)D_0(c)-bD_0(a)D_0(c)\\
&=&(aD_0(b)-bD_0(a))D_0(c).\end{eqnarray*} Then $[aD_0, bD_0]=(aD_0(b)-bD_0(a))D_0\in {\cal D}_1,$
and so 3) holds.  \hfill$\Box$

 \noindent{\bf Theorem 2.4}\quad
{\it Let $\cal A$ be a commutative and associative
algebra satisfying Lemma 2.1, and let  $a$ be an element of $\cal A$.  If $D_0$ satisfies  Lemma 2.3   and $\circ$ satisfies
$$
b\circ c=baD_0(c),  \forall  b, c \in \cal A,\eqno(2.4)$$
then the following statements hold:

$1)$ \quad $({\cal A}, aD_0, \circ)$ is a Novikov algebra.

$2)$ \quad $({\cal A},aD_0, [,])$ is  an adjoining Lie algebra of $({\cal A},
aD_0, \circ)$  and  $[\ ,\ ]$  such that
$$[b, c]=a(bD_0(c)-cD_0(b)),\quad \forall b, c\in \cal A.\eqno(2.5)$$}
{\it Proof}.\quad 1)\, By Lemma 2.3, $aD_0$ is a derivation of the commutative algebra $\cal A$.
 So $({\cal A}, aD_0, \circ)$ is a Novikov algebra by  Remark 1.3 (2).

 2)\quad  $({\cal A}, aD_0, [,])$ is  an adjoining Lie algebra of $({\cal A},
aD_0, \circ)$ by  Remark 1.3 (1). For  $\forall   b, c\in {\cal A},  \exists  a \in \cal A,$ we have
$$ [b,c] =b\circ c-c\circ b= baD_0(c)-caD_0(b) = a (bD_0(c)-cD_0(b))$$ since  $\cal A$ is commutative.
 Hence we obtain the desired result.  \hfill$\Box$

Let $b_0$ be  a unity of $\cal A$. If we set $a=b_0$ in Theorem 2.4, then $
 a_n \circ a_m=a_n b_0D_0(a_m)=a_n (mb_m)=\frac{m}{2}(a_{m+n}+ a_{n-m})$. Similarly,
 we obtain the following corollary:

 \noindent{\bf Corollary 2.5}\quad
 {\it Let $\cal A$ be a commutative and associative
algebra satisfying  Lemma 2.1. Then the following statements hold:
$$\left\{\begin{array}{l}a_n\circ a_m=\frac{m}{2}(a_{n+m}+a_{n-m})\cr
b_n\circ b_m=\frac{m}{2}(a_{n+m}+a_{m-n})\cr a_n\circ
b_m=\frac{m}{2}(b_{n+m}-b_{n-m})\cr b_n\circ
a_m=\frac{m}{2}(b_{n+m}+b_{n-m})\cr\end{array}\right.\eqno(2.6)$$ and
$$\left\{\begin{array}{l}[a_n, a_m]=\frac{1}{2}(m-n)a_{n+m}+\frac{1}{2}(m+n)a_{n-m}\cr
[b_n, b_m]=\frac{1}{2}(m-n)a_{n+m}-\frac{1}{2}(m+n)a_{n-m}\cr
[a_n, b_m]=\frac{1}{2}(m-n)b_{n+m}-\frac{1}{2}(n+m)b_{n-m}\cr
[b_n, a_m]=\frac{1}{2}(m-n)b_{n+m}+\frac{1}{2}(m+n)b_{n-m}.\cr\end{array}\right.\eqno(2.7)$$}

 The following, let $\sinh x=\frac{e^x-e^{-x}}{2}$, $\cosh x=\frac{e^x+e^{-x}}{2}$ the field $\bf F$ be assumed $\bf R$ or $\bf C$. We will construct Novikov algebras
 over the linear space which is generated by  $\sinh x$ and $\cosh x$.

     First, let $\cal T$ be a linear space generated by $\{\sinh mx, \cosh nx|m,n\in{\bf N}\}$
     over   $\bf F$.

\noindent{\bf Lemma 2.6}\quad
 {\it $\cal T$ satisfying the above product is a commutative  associative algebra.}

{\it Proof.}\,
Since the  above product   is  commutative and
 associative, we only need that $\cal T$ is closed for the product.
In fact,
$$
\left\{\begin{array}{l}\sinh mx\sinh nx=\displaystyle\frac{1}{2}[\cosh (m+n)x-\cosh (m-n)x]\cr \cr
 \cosh mx\cosh nx=\displaystyle\frac{1}{2}[\cosh (m+n)x+\cosh (m-n)x]\cr  \cr
\sinh mx\cosh nx=\displaystyle\frac{1}{2}[\sinh (m+n)x+\sinh
(m-n)x].\cr\end{array}\right.\eqno(2.8)$$ So $\cal T$ is a
commutative and associative algebra.  \hfill$\Box$

   \noindent{\bf Lemma 2.7}\quad {\it Let $\cal T$ be a linear space generated by $\{\sinh mx, \cosh
nx|m, n\in{\bf N}\}$ over $\bf F$, then $\{1, \sinh mx, \cosh nx|m,
n\in {\bf N_0}\}$ is a basis of  $\cal T$.}

 {\it Proof.}\quad For $\forall n\in {\bf N}_0$, suppose that there
 are $c_0,
a_i, b_j\in {\bf F}, i, j\in {\bf N_0}$ such that
$$
c_0+a_1\sinh x+b_1\cosh x+\cdots+a_n\sinh nx+b_n\cosh
nx=0.\eqno(2.9)$$ We take derivative for (2.9) such that its
derivative order is $2k-1$ $(k\in{\bf N}_0)$, and put $x=0$. Then we
have
$$a_1+2^{2k-1}a_2+\cdots+n^{2k-1}a_n=0.$$
Let $k=1, 2, \cdots, n$, then we obtain the following system of $n$
linear equations:  $$ \left\{\begin{array}{l}
a_1+2a_2+\cdots+na_n=0\\
a_1+2^3a_2+\cdots+n^3a_n=0\\
\cdots\cdots\cdots\cdots\cdots\cdots\cdots\cdots\cdots\\
a_1+2^{2n-1}a_2+\cdots+n^{2n-1}a_n=0.
\end{array}\right.\eqno(2.10)$$
If $a_1, \cdots, a_n$ are seen to be unknown,  then the coefficient
matrix of (2.10) is the Vandermonde matrix whose determinant is not
$0$, so $a_i=0, i=1, \cdots, n$.

We take derivative for (2.9) such that its derivative order is $2k$
$(k\in{\bf N}_0)$, and put $x=0$. Then we have
$$b_1+2^{2k}b_2+\cdots+n^{2k}b_n=0.$$
Let $k=1, 2, \cdots, n$, then we obtain the following system of $n$
linear equations:$$ \left\{\begin{array}{l}
b_1+2^2b_2+\cdots+n^2b_n=0\\
b_1+2^4b_2+\cdots+n^4b_n=0\\
\cdots\cdots\cdots\cdots\cdots\cdots\cdots\cdots\cdots\\
b_1+2^{2n}b_2+\cdots+n^{2n}b_n=0.
\end{array}\right.\eqno(2.11)$$
If $b_1, \cdots, b_n$ are seen to be unknown,  then the coefficient
matrix of (2.11) is the Vandermonde matrix whose determinant is not
$0$, so $b_i=0, i=1, \cdots, n$. Since for $\forall i\in {\bf N}_0$,
$a_i=0$ and $b_i=0$ satisfy (2.9), we have $c_0=0$. Hence $\{1$,
$\sinh x$, $\cosh x$, $\cdots$, $\sinh nx$, $\cosh nx\}$ are
linearly independent for $\forall n\in {\bf N}_0$, then $\{1, \sinh
nx, \cosh mx|n, m\in {\bf N}_0\}$ are linearly independent and so
they form a basis of $\cal T$ as desired. \hfill$\Box$

 \noindent{\bf Theorem 2.8}\quad
 {\it Let ${\cal A}_1$, ${\cal A}_2$  be
commutative and associative algebras over  $\bf F$. If $\varphi$:
${\cal A}_1\longrightarrow {\cal A}_2$ is an isomorphism and $D_1\in
{\rm Der}{\cal A}_1$, then the following statements hold:

$1)$\quad $D_2: =\varphi D_1\varphi^{-1}\in {\rm Der}{\cal A}_2$.

$2)$\quad $\varphi$: $({\cal A}_1, D_1,\circ ) \longrightarrow({\cal
A}_2, D_2, \circ)$ is also an isomorphism of Novikov algebras.}

{\it Proof.}\quad
1)\quad For $\forall  a, b\in{\cal A}_1$, we have
\begin{eqnarray*}&&\ \ (\varphi D_1\varphi^{-1})(\varphi(a)\varphi(b))=
(\varphi D_1\varphi^{-1})(\varphi(ab))\\
&&=\varphi D_1(ab)=\varphi (D_1(a)b+aD_1(b))=\varphi(D_1(a))\varphi(b)+\varphi(a)\varphi(D_1(b))\\
&&=(\varphi
D_1\varphi^{-1})(\varphi(a))\varphi(b)+\varphi(a)(\varphi
D_1\varphi^{-1})(\varphi(b)).
\end{eqnarray*}
So 1)  holds.

2)\quad For $ \forall a, b\in{\cal A}_1$, we have
\begin{eqnarray*}&&\ \ \varphi(a\circ b)=\varphi(aD_1(b))=\varphi(a)\varphi(D_1(b))\\
&&=\varphi(a)(\varphi D_1\varphi^{-1})(\varphi(b))=\varphi(a)D_2(\varphi(b))\\
&&=\varphi(a)\circ \varphi(b).\end{eqnarray*} So 2)  holds.  \hfill$\Box$

 \noindent{\bf Theorem 2.9}\quad
{\it Let ${\cal A}$ be  a commutative and associative algebra over
$\bf F$  satisfying  Lemma 2.1, $D_0$ be its derivation satisfying
(2.2) and ${\cal T}$ be a commutative and associative algebra over
$\bf F$ satisfying Lemmas 2.6 and 2.7. If  $\varphi: {\cal A}
\longrightarrow {\cal T}$
  satisfies
$$
\varphi(b_m)=\cosh mx,\ m=0, 1, 2, \cdots, \quad \varphi(a_n)=\sinh
nx,\ n=1, 2, \cdots,\eqno(2.12)$$ then the following statements
hold:

$1)$\quad  $\varphi$ is an  isomorphism of commutative and associative algebras.

$2)$\quad $\varphi D_0\varphi^{-1}=\frac{\rm d}{{\rm d}x}$.

$3)$\quad $\varphi: ({\cal A}, aD_0, \circ)\longrightarrow ({\cal
T}, \varphi(a)\frac{\rm d}{{\rm d}x}, \circ)$ is an isomorphism of
Novikov algebras.}

 {\it Proof.}
It is clear by  Lemma 2.7,  (2.1) and (2.8).

2) By (2.2) and (2.12), we have
\begin{eqnarray*}&&\ \  \varphi D_0 \varphi^{-1}(\sinh nx)= \varphi D_0(a_n)\\
&&= \varphi(nb_n)=n\cosh nx\\
 &&=\frac{{\rm d}\sinh nx}{{\rm d}x},\\
&&\ \  \varphi D_0 \varphi^{-1}(\cosh nx)= \varphi D_0(b_n)\\
&&= \varphi(na_n)=n\sinh nx\\
 &&=\frac{{\rm d}\cosh nx}{{\rm d}x}.\end{eqnarray*}
 So 2)  holds.

 3)\quad It is clear that $\varphi ( aD_0)\varphi^{-1} = \varphi(a)d/dx $.
 By (2.12) and (2.2), we have
 \begin{eqnarray*}&&\ \
       \varphi( aD_0)\varphi^{-1} ( \sinh nx)
         =\varphi (aD_0)( a_n)\\
        && =\varphi(a D_0(a_n))
        =\varphi( anb_n) \\
        &&=\varphi(a)\varphi(nb_n)
      = \varphi(a)n \cosh nx\\
     &&=\varphi(a)d(\sinh nx)/dx. \end{eqnarray*}
  Similarly, we have $\varphi ( aD_0)\varphi^{-1} ( \cosh nx) = \varphi(a)d(\cosh nx)/dx$.
 So $\varphi( aD_0)\varphi^{-1}= \varphi(a)d/dx.$

    By  Theorems 2.4, 2.8 and Remark 1.3 (2), we have
  \begin{eqnarray*}&&\ \
\varphi( b\circ c)=\varphi(baD_0(c)) \\&&=\varphi(b)\varphi(aD_0(c)) \\
&&=\varphi(b)[ \varphi(aD_0)\varphi^{-1}( \varphi(c)) ]\\
 &&=\varphi(b) \varphi(a)d/dx ( \varphi(c))\\
 &&=\varphi(b)\circ \varphi(c), \forall b,c\in {\cal A}. \end{eqnarray*}
So $\varphi: ({\cal A}_0, aD_0, \circ)\longrightarrow ({\cal
T}, \varphi(a)\frac{\rm d}{{\rm d}x}, \circ)$ is an isomorphism of
Novikov algebras.    \hfill$\Box$

\end{document}